\documentclass{article}

\usepackage{amssymb,amsfonts,enumerate,version,amsmath,theorem,a4wide}
\usepackage[T1]{fontenc}
\usepackage[math]{iwona}
\usepackage{hyperref}

\theoremstyle{change}
\newtheorem{Theorem}{Theorem}[section]

\newcommand\caps{\scshape}

\renewcommand\L{\mathcal{L}}

\newcommand\half{\frac12}

\newcommand\inv{^{-1}}

\def\pointerstart{\color{red}\ \pointersize \Pointinghand\ \ }
\def\pointersize{\Large}
\def\ppointer#1\par{\pointerstart
{\fbox{
    \raisebox{2ex}{\parbox[t]{9cm}{\color{blue}\pointersize #1}}}}
}

\newcommand\norm[1]{\Vert #1 \Vert}

\newcommand\summ[3]{\displaystyle\sum_{#1 = #2}^{#3}}

\newcommand\ip[2]{\langle #1, #2\rangle}

\newcommand\BC[2]{\binom {#1}{#2}}

\newcommand\tp[3]{\{#1,#2,#3\}}
\newcommand\D[2]{#1 {\scriptsize\square} #2}

\newtheorem{Ex}[Theorem]{Example}
\newtheorem{Def}[Theorem]{Definition}
\newtheorem{Lemma}[Theorem]{Lemma}
\newtheorem{Proposition}[Theorem]{Proposition}
\newtheorem{Corollary}[Theorem]{Corollary}
\newtheorem{Conjecture}[Theorem]{Conjecture}
\newtheorem{kNote}[Theorem]{Note}
\newtheorem{Nota}[Theorem]{Notation}
\newtheorem{Exer}[Theorem]{Exercise}
\newtheorem{Rem}[Theorem]{Remark}
\newtheorem{kNotes}[Theorem]{Notes}

\newenvironment{Example}{\begin{Ex} \rm}{\end{Ex}}
\newenvironment{Definition}{\begin{Def}\rm}{\end{Def}}

   \def\exname{ex}

\newenvironment{Proof}{\noindent\textbf{Proof.}}{\hfill$\square$\medskip}

\renewcommand\j{\circ}
\newcommand\q[1]{Q_{#1}}

\title{Local Derivations on Jordan Triples}
\author{Michael Mackey\footnote{UCD School of Mathematical Sciences \hfill email: michael.mackey@ucd.ie}}

\begin{document}

\maketitle

\begin{abstract}
  R.V. Kadison defined the notion of local derivation on an
  algebra and proved that every continuous local derivation on a von Neumann
  algebra is a derivation~\cite{MR1051316}.  We provide the analogous result in the
  setting of Jordan triples.
\end{abstract}

\makeatletter
  \def\footnotestar{\xdef\@thefnmark{}\@footnotetext}
\makeatother

\footnotestar{AMS Subject Classification (2000): 17C65, 47Cxx}

\section{Introduction}

R.V.~Kadison gave the definition of local derivation on an associative
algebra and proved the fundamental result that
every continuous local derivation on a von Neumann algebra is a derivation.  In
the intervening period, a substantial body of literature has been
built up on the topic of local derivations and local automorphisms.
As example, B.E.~Johnson~\cite{MR1783788} extended
Kadison's theorem to derivations on  arbitrary C*-algebras and also
showed that the 
continuity of such local derivations is automatic. 

Jordan algebras, Jordan Banach algebras and JC*-algebras are a widely studied generalisation of their associative
counterparts.  These structures are further subsumed into triple
analogues: Jordan triples, Jordan Banach triples and JC*- and
JB*-triples, which can be loosely interpreted as ``rectangular''
versions of their binary or ``square'' forebears. Definitions follow
below.   It is our intention in this note to extend the main result
of~\cite{MR1051316}
to the the setting of Jordan triples.  Let us make some initial limitations
to our study.  Firstly, while the original result of Kadison (and that
of Johnson) related
to module-valued maps, we deal only with self-maps on the Jordan
triple.  Modules are not commonly considered in Jordan triple theory,
one reason (other than the algebraic difficulty) being that,
frequently,
 a particular module over a (Jordan) algebra may be, 
itself, a (Jordan) triple and so passing to the triple setting
obviates the need to consider modules.  As examples, the right module
$\L(H,K)$ over the C*-algebra $\L(H)$ is a JB*-triple, and every Hilbert C*-module is a
JB*-triple~\cite{MR1997703}.

Let us point out that, while Jordan triple
structures generalise their binary associative counterparts, specific
properties of elements or mappings may not.  For example, the triple
analogue of a (binary) idempotent is known as a tripotent, but an
idempotent of an algebra may not be a tripotent when that algebra is
viewed as a triple.  Another example, most pertinent to us, is that a
derivation on an algebra may not be a triple derivation when the
algebra is considered as a Jordan triple.  Indeed, we
examine this aspect more closely in Section~\ref{sec:ders}.  Thus,
while our
main result is, in spirit, a generalisation of that of Kadison, it
also provides something new in the category of von Neumann algebras.

\section{Background and Terminology}
\label{sec:bg}

For the reader unfamiliar with the notions of Jordan algebra and
Jordan *-algebra we refer to \cite{HOS} or \cite{McCrimmon_taste}.
 
\begin{Definition}
  \begin{enumerate}
    \item 
Let $Z$ be a complex vector space.  A \textsc{triple product} on $Z$
is a real tri-linear map $\tp\cdot\cdot\cdot:Z^3\to Z$, $(x,y,z)\mapsto \tp xyz$, which is 
      \begin{enumerate}
      \item complex linear in the outer variables $x$ and $z$,
        \item symmetric in the outer variables, that is $\tp xyz =\tp
          zyx$, and
          \item complex anti-linear in the inner variable $y$.
      \end{enumerate}
\item A \textsc{Jordan triple} is a complex vector space with triple
  product which satisfies the Jordan
      triple identity
      \begin{equation}
        \label{eq:JTI}
        \tp ab{\tp xyz} = \tp {\tp abx}yz - \tp x{\tp bay}z + \tp
        xy{\tp abz}.
      \end{equation}
    \end{enumerate}
  \end{Definition}
We say that a triple product is \textsc{non-degenerate} if $\tp xxx
=0$ only when $x=0$.
Every algebra becomes a Jordan algebra under the Jordan product
defined by $x\circ y = \half(x.y+y.x)$ while every Jordan *-algebra is a
Jordan triple via the triple product 
\begin{equation}\label{eq:jtt}
\tp xyz = (x\circ y^*)\circ z +
(z\circ y^*)\circ x - (x\circ z)\circ y^*.
\end{equation}
Combining these two
facts, one sees that every *-algebra is a Jordan triple via the triple
product
\begin{equation}
\tp xyz = \half(xy^*z+ zy^*x).\label{eq:att}
\end{equation}

By \textsc{derivation} on an algebra, we mean a linear map $d$
satisfying $d(x.y)=d(x).y + x.d(y)$.  A \textsc{Jordan derivation} on
an (associative) algebra is a linear map $d$ for which
$d(x^2)=x.d(x)+d(x).x$ (or equivalently and in terms of the Jordan
product, $d(x\circ y)= x\circ d(y)+ d(x)\circ y$).  Also in the
literature one finds the notion of \textsc{Jordan *-derivation} on an
algebra: a \emph{real} linear map $d$ for which
$d(x^2)=x.d(x)+d(x).x^*$.  Respective examples of these are $x\mapsto
ax-xa$ and $x\mapsto ax-x^*a$.  A derivation of a Jordan algebra is
again a linear map $d$ with $d(x\circ y)= x\circ d(y)+d(x)\circ y$.
Note that a Jordan derivation on an algebra is a derivation on the
associated Jordan algebra.

The ensuing concept for a Jordan
triple is natural.

\begin{Definition}
A \textsc{(Jordan) triple derivation} is a linear map $d$ on a Jordan triple
  satisfying \[d\tp xyz = \tp{dx}yz + \tp x{dy}z +\tp xy{dz}.\]
\end{Definition}

The importance of the notion of triple derivation in the framework of
Jordan triples is apparent upon noting that the Jordan triple identity \eqref{eq:JTI} can be equivalently formulated
thus: for all $a$ and $b$, the map $\D ab - \D ba$ is a triple
derivation, where $\D ab$ denotes the linear operator $x\mapsto\tp abx$.  Indeed, a further reformulation in our complex setting
is that for all $a$, $i\D aa$ is a triple derivation.  Derivations of
this form are known as \emph{inner} derivations.

By a Jordan triple derivation on a Jordan *-algebra or an associative
*-algebra, we mean with respect to the triple product as given
by~\eqref{eq:jtt} and \eqref{eq:att} respectively.
Jordan triple derivations form a \emph{real} linear space and are closed
under the Lie bracket- that is if $d_1$ and $d_2$ are Jordan triple
derivations then so is $d_1d_2-d_2d_1$.  We refer to \cite{MR1067482} and
\cite{MR1875137} for greater detail converning the structure of triple
derivations on JB*-triples, these forming a quite specialised class of
non-degenerate Jordan triple, which we now introduce.

\begin{Definition}
 A \textsc{JB*-triple} is a complex Banach space and a Jordan triple on
     which the triple product is jointly continuous and satisfies for
     every element $x$:
     \begin{enumerate}[(i)]
     \item $\sigma(\D xx) \ge 0$,
       \item $\exp(i\D xx)$ is a triple automorphism and a
         surjective linear isometry,
         \item $\norm {\D xx(x)} = \norm x^3$. 
     \end{enumerate}
\end{Definition}

The class of JB*-triples includes all C*-algebras (via $\tp xyz=\half
(xy^*z+zy^*x)$) and also Hilbert space (via $\tp xyz=\half(\ip xyz
+\tp zyx$) and $\mathcal{L}(H,K)$ where $H$ and $K$ are Hilbert spaces.  The category is more than a convenient envelope for
better known structures however- it stands independently as a complex analytic
category because of the following theorem proven by
W.~Kaup~\cite{Kaup_RMT} following a body of work that can be traced
back\footnote{We mention here the names of Jordan, von Neumann,
  Wigner, Tits, Kantor, Koecher, Loos, and Harris.} to Elie Cartan's classification of Hermitian symmetric spaces~\cite{Cartan}.  
\begin{Theorem}
  A Banach space is a JB*-triple if, and only if, its open unit ball
  has a transitive group of biholomorphic mappings.
\end{Theorem}
Let us mention some of the principle tools and relevant facts when
working with a JB*-triple $Z$. For all $x,y,z \in Z$, $\norm{\tp
  xyz}\le \norm x \norm y \norm z$ (the proof of which~\cite{FR_GN} does not follow easily
from the definitions).
The linear operator $B(x,y)\in\mathcal{L}(Z)$ defined by \[B(x,y)= I - 2\D
xy + \q x\q y\] where $\q x(z) = \tp xzx$ occurs frequently and is
known as the Bergman operator.  On a C*-algebra, the Bergman operator
reduces to $B(x,y) z = (1-xy^*)z(1-y^*x)$.  Derivations of JB*-triples
are automatically bounded \cite{MR1067482}.

An element $e\in Z$ for
which $\tp eee=e$ is called a tripotent and a
non-zero tripotent has norm one.  For example, a tripotent of a C*-algebra
is an element $v$ satisfying $v=vv^*v$, that is, a partial isometry.
Each tripotent induces a splitting of $Z$, called the Peirce
decomposition, into $Z=Z_1 \oplus Z_\half \oplus Z_0$ where $Z_k$ is
the $k$-eigenspace of $\D ee$, with mutually orthogonal projections
$P_k$ onto the subspaces $Z_k$,
\begin{align*}
        P_1 &= \q e\q e, \\
        P_\half &= 2\D ee -2\q e\q e,\\
        P_0 &= B(e,e),
\end{align*}
satisfying $P_1+P_\half + P_0 = I$.  Where the need arises, we write
$P^e_j$ rather that $P_j$ to highlight the tripotent in question.
With respect to this decomposition, the triple product behaves as
follows:
\[ \tp{Z_i}{Z_j}{Z_k} \subset Z_{i-j+k}\]
where $i,j,k\in\{0,\half,1\}$ and $Z_p=\{0\}$ when
$p\notin\{0,\half,1\}$.  In addition, $\D{Z_0}{Z_1}=0=\D{Z_1}{Z_0}$.  The
tripotent $e$ is called {\caps maximal} if $Z_0=\{0\}$ and this is the
case precisely when $e$ is an extreme point of the unit ball of $Z$
\cite{KaupUpmeier_Siegel}.  We
point out that since the triple product is continuous, the set of
tripotents forms a closed subset of $Z$.  Two tripotents are said to
be orthogonal if $\D ef=0$ which is equivalent to saying $\D fe=0$ or
that $\tp eef=0$.  Note that the sum of two orthogonal tripotents is a tripotent.

The bidual of a JB*-triple is also a JB*-triple (\cite{Dineen_CHVF}) and
any JB*-triple with a (necessarily unique) predual is called a JBW*-triple
(cf.~\cite{BT}).  While a JB*-triple may not have any
tripotents, a JBW*-triple has an abundance.  For example, the unit ball of
a JBW*-triple is the convex hull of its (maximal) tripotents.
In addition, we have~\cite[Lemma 3.11]{MR929400}:

\begin{Proposition}\label{prop:dens}
  The set of tripotents is norm total in a JBW*-triple.  More
  precisely, each element in the JBW*-triple can be approximated in
  norm by  a finite linear combination of mutually orthogonal tripotents.
\end{Proposition}

There is, as in the algebra setting, a close link between Jordan
triple derivations
and Jordan triple automorphisms. 

\begin{Lemma}\label{lem:expder}
  Suppose $d$ is a derivation on a JB*-triple.  Then $\exp d = \summ
  n0\infty d^n/n!$ preserves the Jordan triple product.  Conversely,
  if $\exp td$ is a triple homomorphism for all $t>0$.  Then $d$
  is a derivation.
\end{Lemma}




There is another quite particular reason why, in the context of
JB*-triples, derivations are of special interest.  As proven by
Kaup~\cite{Kaup_RMT}, the Jordan triple
automorphisms, that is the bijective bounded linear maps which preserve the triple
product, coincide precisely with the surjective linear isometries.
That is, for $T\in \mathsf{GL}(Z)$, $\norm{Tx}=\norm x$ for all $x\in
Z$ if and only if $T\tp xyz = \tp{Tx}{Ty}{Tz}$ for all $x,y,z\in
Z$. (When the surjectivity requirement is dropped, the situation is
rather more complicated however, cf.~\cite{MR2190348}.)
In light of Lemma~\ref{lem:expder} therefore, it is valid to think of a
triple derivation as the ``infinitesimal'' form of an isometry.  Also
worthy of note here is the fact that, on each of the irreducible JB*-triples
known as the Cartan Factors (see e.g.~\cite{FR_GN}), the group of inner automorphisms (that is,
those in the group generated by exponentials of inner derivations)
acts 
transitively on the manifold of tripotents of a given
rank~\cite{MR2491603}.  (The rank 
of a tripotent $e$ is the maximum number of mutually orthogonal tripotents
whose sum equals $e$.)

\section{Jordan derivations and Jordan triple derivations}
\label{sec:ders}

We are interested in clarifying how Jordan triple derivations relate to
Jordan derivations.   In particular, what algebraic conditions on a
linear map $d:A\to A$ (where $A$ is a *-algebra) are equivalent
to $d$ being a Jordan triple derivation.

The reason behind this question is that neither Jordan derivations,
*-preserving Jordan derivations nor Jordan *-derivations on an
associative *-algebra $A$ produce triple derivations on the associated
Jordan triple.  This contrasts with, say, Jordan triple automorphisms,
which certainly generalise associative *-automorphisms.  In
particular, notice that if one assumes the *-algebra $A$ has an
identity $1$ (as we will throughout since any derivation on the non-unital
algebra $A$ extends to the unitisation by linearity and $d(1):=0$)
then a Jordan derivation, or Jordan *-derivation on $A$ sends the
identity to $0$.  However the Jordan triple derivation $i\D aa$ is generally
non-zero, as one can see in the definition of a JB*-triple.

Let us begin with a simple observation.
\begin{Lemma}\label{lem:a} Let $\delta$ be a Jordan triple derivation on a
(unital) Jordan *-algebra. Then
\begin{enumerate}[(a)]\item 
  $  \delta(a\j b) = \delta a\j b + a\j \delta b + \tp a{\delta1}b,$
  \item $\delta(b^*) = 2\delta1\j b^* + (\delta b)^*$.
\end{enumerate}

\end{Lemma}

\begin{Proof} For the first part, write $a\circ b=\tp a1b$ and
  apply $\delta$.  For the second, use $b^*=\tp 1b1$.
\end{Proof}

We use $L_x$ to denote the multiplication operator on a Jordan
algebra, that is $L_x(y)=x\circ y$.
\begin{Lemma}\label{lem:b}
  Let $A$ be a unital Jordan *-algebra.  If $x=-x^*$ then $L_{x}$ is a Jordan triple derivation.  Further, the
  converse holds if the
  triple product is non-degenerate.
\end{Lemma}

\begin{Proof}
  \begin{align}
    L_{x}(\tp abc) &= x \j \tp abc \notag\\
                        &= \tp {x}1{\tp abc} \notag\\
                        &= \tp{\tp{x}1a}bc - \tp a{\tp
                          1{x}b}c + \tp ab{\tp{x}1c} \notag\\
                        &= \tp {L_{x}a}bc + \tp ab{L_{x}c}
                        - \tp a{\tp 1{x}b}c.\label{eq:1}
  \end{align}
Since  $\tp 1xb = x^*\circ b = -x\circ b =
-L_x b$ we have that $L_x$ is a Jordan triple
derivation as required.

Towards the converse, we see from \eqref{eq:1} that if $L_x$ is a
triple derivation then for all $a,b$ and $c$, $\tp a{\tp x1b}c = -\tp
a{\tp 1xb}c$.  The non-degeneracy of the triple product implies that
$\D x1 = -\D 1x$ and in particular that $x^*=\tp 1x1 = -\tp x11 = -x$.
\end{Proof}

The following corollary of this fact appears in \cite[Lemma 1]{MR1875137}.
\begin{Corollary}\label{cor:a}
  If $\delta$ is a triple derivation on a Jordan *-algebra then
  so is $L_{\delta1}$.
\end{Corollary}

\begin{Proof}
  This follows from the fact that $\delta1=\delta\tp 111 =
  2\tp{\delta1}11 + \tp 1{\delta1}1 = 2\delta1 + (\delta1)^*$.  Thus
  $(\delta1)^*=-\delta1$ and we can apply the previous result.
\end{Proof}

\begin{Theorem}\label{thm:char}
  A linear map $\delta$ on a Jordan *-algebra is a Jordan triple
  derivation if, and only if,
  \begin{enumerate}[(i)]
  \item $\delta(a\j b) = \delta a\j b + a\j \delta b + \tp
    a{\delta1}b$ for all $a$ and $b$, and
  \item $\delta(b^*) = 2\delta1\j b^* +(\delta b)^*$ for all $b$.
  \end{enumerate}
\end{Theorem}

\begin{Proof}
  The necessity of the two conditions is provided by~\ref{lem:a}.
  So suppose $\delta$ satisfies the conditions above.  From
  \eqref{eq:jtt} we have
  \[ \delta\tp abc = \delta((a\j b^*)\j c) + \delta((c\j b^*)\j a) -
  \delta((a\j c)\j b^*).\]
  By use of (i), we expand as follows.
  \begin{align*}
    \delta((a\j b^*)\j c) &= \delta(a\j b^*)\j c + (a\j b^*)\j \delta c
    + \tp {a\j b^*}{\delta1}c \\
      &= [\delta a\j b^* + a\j\delta(b^*) + \tp a{\delta1}{b^*}]\j c +
      (a\j b^*)\j \delta c + \tp{a\j b^*}{\delta1}c \\
      &= (\delta a\j b^*)\j c + (a\j \delta(b^*))\j c + +\tp
      a{\delta1}{b^*}\j c +(a\j b^*)\j
      \delta c + \tp{a\j b^*}{\delta1}c 
  \end{align*}
  Similarly we have
  \begin{align*}
    \delta((c\j b^*)\j a) &= (\delta c\j b^*)\j a + (c\j \delta(b^*))\j a + +\tp
      c{\delta1}{b^*}\j a +(c\j b^*)\j
      \delta a + \tp{c\j b^*}{\delta1}a\\
      \intertext{and}
    \delta((c\j a)\j b^*) &=
       (\delta c\j a)\j b^* + (c\j \delta a)\j b^* + +\tp
      c{\delta1}{a}\j b^* +(c\j a)\j
      \delta b^* + \tp{c\j a}{\delta1}{b^*}.
  \end{align*}
  Thus
  \begin{align*}
\delta\tp abc &= (\delta a\j b^*)\j c + (a\j \delta(b^*))\j c + \tp
      a{\delta1}{b^*}\j c +(a\j b^*)\j
      \delta c + \tp{a\j b^*}{\delta1}c \\ &\quad + 
         (\delta c\j b^*)\j a + (c\j \delta(b^*))\j a + \tp
      c{\delta1}{b^*}\j a +(c\j b^*)\j
      \delta a + \tp{c\j b^*}{\delta1}a \\ &\quad - 
         [(\delta c\j a)\j b^* + (c\j \delta a)\j b^* + +\tp
      c{\delta1}{a}\j b^* +(c\j a)\j
      \delta b^* + \tp{c\j a}{\delta1}{b^*}]\\ \\
      &= \tp {\delta a}bc + \tp a{(\delta (b^*))^*}c + \tp ab{\delta
        c} \\
      &\qquad + \tp a{\delta1}{b^*}\j c + \tp {a\j
          b^*}{\delta1}c + \tp c{\delta1}b^* \j a + \tp {c\j
          b^*}{\delta1}a \\ &\qquad- \tp c{\delta1}a \j b^* - \tp{c\j
          a}{\delta1}{b^*}. 
  \end{align*}
From condition (ii), 
$(\delta(b^*))^* = \delta b + 2b\j (\delta1)^*$ and (on taking $b=1$)
$(\delta1)^*=-\delta1$.  Therefore $(\delta(b^*))^* = \delta b - 2b\j
\delta1$ which we use with the above to write
\begin{align}
  \delta\tp abc &= \tp {\delta a}bc + \tp a{\delta b}c + \tp ab{\delta
    c} \notag\\&\quad -2\tp a{b\j \delta1}c +  \biggl[ 
       \tp a{\delta1}{b^*}\j c + \tp {a\j
          b^*}{\delta1}c + \tp c{\delta1}b^* \j a + \tp {c\j
          b^*}{\delta1}a \notag \\ &\qquad- \tp c{\delta1}a \j b^* - \tp{c\j
          a}{\delta1}{b^*} \biggr].\label{eq:3}
\end{align}
Consider the square-bracketed terms in this expression.  Expanding via
the algebra product we have
\begin{align*}
  \tp a{\delta1}{b^*}\j c &=((a\j (\delta1)^*)\j {b^*} + ({b^*}\j {(\delta1)^*})\j a
  - (a\j {b^*})\j {(\delta1)^*})\j c \\
\tp{a\j {b^*}}{\delta1}c &= ((a\j {b^*} )\j {(\delta1)^*})\j c + (c\j {(\delta1)^*})\j
(a\j {b^*}) - ((a\j b^*)\j c)\j {(\delta1)^*} \\
\tp c{\delta1}{b^*} \j a &= ((c\j (\delta1)^*)\j {b^*} + ({b^*}\j {(\delta1)^*})\j c
  - (c\j {b^*})\j {(\delta1)^*})\j a \\
\tp{c\j {b^*}}{\delta1}a &= ((c\j b^*)\j {(\delta1)^*})\j a + (a\j
{(\delta1)^*})\j(c\j {b^*}) - ((c\j b^*)\j a)\j {(\delta1)^*} \\
-\tp c{\delta1}a\j {b^*} &= -b^*\j (  (c\j{(\delta1)^*})\j a + (a\j {(\delta1)^*})\j
c  - (a\j c)\j {(\delta1)^*}) \\
-\tp{c\j a}{\delta1}{b^*} &= -((c\j a)\j {(\delta1)^*})\j b^* - (b^*\j {(\delta1)^*})\j
(c\j a) + ((c\j a)\j b^*)\j {(\delta1)^*}
\end{align*}
Summing these, we see the
square-bracketed terms in~\eqref{eq:3} can be written as
\[ \tp{a\j(\delta1)^*}bc + \tp ab{c\j(\delta1)^*} + \tp a{b\j\delta1}c
- \tp abc\j(\delta1)^*.\]  Replacing $(\delta1)^*$ by $-\delta1$, we
substitute into~\eqref{eq:3} to get
\begin{align*} 
\delta\tp abc &= \tp{\delta a}bc + \tp a{\delta b}c + \tp ab{\delta c}
\\ &\qquad -2\tp a{b\j \delta1}c  - \tp{a\j\delta1}bc -\tp ab{c\j\delta1} + \tp
a{b\j\delta1}c + \tp abc\j\delta1\\
&= \tp{\delta a}bc + \tp a{\delta b}c + \tp ab{\delta c} \\ &\qquad- \tp{a\j\delta1}bc -\tp ab{c\j\delta1} - \tp
a{b\j\delta1}c + \tp abc\j\delta1.
\end{align*}
Again since $(\delta1)^*=-\delta1$, Lemma~\ref{lem:b} gives that $L_{\delta1}$ is a
Jordan triple derivation and so the above reduces to  
\[\delta\tp abc = \tp{\delta a}bc + \tp a{\delta b}c + \tp ab{\delta c}\]
as required.
\end{Proof}


\begin{Example}
  Let $A$ be an associative *-algebra with the usual Jordan binary and
  triple products.  For $\delta=M_{a,b}$ defined by \[ M_{a,b}(x)=xa +
  bx\] 
  condition (i) is equivalent to $(a+b)^*= -(a+b)$
  while condition (ii) is equivalent to $x(a+a^*)+(b+b^*)x = 0$ for
  all $x\in A$.
\end{Example}

\section{Derivation Pairs}

Recall that a Jordan triple isomorphism $\lambda$, is an invertible linear map which
preserves the triple product, $\lambda\tp xyz = \tp {\lambda
  x}{\lambda y}{\lambda z}$.  A more general concept is that of
structure map, which is actually a pair of invertible linear maps $(S,T)$ which
satisfy $S\tp x{Ty}z = \tp {Sx}y{Sz}$ and $T\tp x{Sy}z =
\tp{Tx}y{Tz}$, or equivalently, 
\[ S\tp xyz =\tp {Sx}{T\inv y}{Sz} \quad\mathrm{and}\quad T\tp xyz =
\tp {Tx}{S\inv y}{Tx}.\]
Structure maps can be used to define homotopes of Jordan structures \cite{Mackey_homotopes}.   A structure map $(S, T)$ is  a Jordan triple isomorphism when
$S=T\inv$.  Lemma~\ref{lem:expder} prompts us to the following
definition.

\begin{Definition}
  A \textsc{derivation pair} on a Jordan triple is a pair of linear
  maps $D=(d_+,d_-)$ which satisfy
  \begin{gather*}
    d_+\tp xyz = \tp{d_+x}yz+\tp x{d_- y}z + \tp xy{d_+z}\\
    d_-\tp xyz = \tp{d_- x}yz+\tp x{d_+ y}z + \tp xy{d_- z}
  \end{gather*}
for all $x, y $ and $z$.
\end{Definition}

For example, the Jordan triple identity states that $(\D xy, -\D yx)$
is a derivation pair, and that $(i\D xx, i\D xx)$ is a derivation
pair.  Clearly $d$ is a derivation if, and only if, $(d,d)$ is a
derivation pair.  As long as the triple under question is
non-degenerate, $d_+$ and $d_-$ uniquely determine one another when
$(d_+,d_-)$ is a derivation pair.  If $(d^1_+, d^1_-)$ and $(d^2_+,
d^2_-)$ are derivation pairs then so is
$([d^1_+,d^2_+],[d^1_+,d^2_-])$.   Iterative action of the derivation
pair $(d_+,d_-)$ follows expected rules:
\begin{gather*}
  d_{\pm}^n\tp xyz = \summ k0n\summ l0{n-k} \BC nk\BC{n-k}l\tp
  {d_{\pm}^kx}{d_{\mp}^ly}{d_\pm^{n-k-l}z}.
\end{gather*}
Notational modifications of
Lemma~\ref{lem:expder} guarantee the following fact.

\begin{Lemma}
  If $(d_+,d_-)$ is a derivation pair then $(\exp d_+, (\exp d_-)\inv)$ is
  a structure map.
\end{Lemma}

In fact, all the results from Section~2 can be stated in terms of derivation
pairs.

\begin{Lemma}\label{lem:apair} Let $\delta=(\delta_+,\delta_-)$ be a
  Jordan triple derivation pair on a
Jordan *-algebra. Then
\begin{enumerate}[(a)]\item 
  $  \delta_\pm(a\j b) = \delta_\pm a\j b + a\j \delta b + \tp a{\delta_\mp1}b,$
  \item $\delta_\pm(b^*) = 2\delta_\pm1\j b^* + (\delta_\mp b)^*$.
\end{enumerate}
\end{Lemma}

\begin{Lemma}\label{lem:bpair}
  Let $A$ be a unital Jordan *-algebra.  For any $x\in A$, the pair
  $(L_x, L_{-x^*})$ is a Jordan triple derivation pair.  Further, if the
  triple product is non-degenerate then $(L_x, L_y)$ is a derivation
  pair only if $y=-x^*$.
\end{Lemma}

\begin{Corollary}\label{cor:apair}
  If $(\delta_+,\delta_-)$ is a triple derivation on a Jordan *-algebra then
  so is $(L_{\delta_+1}, L_{\delta_-1})$.
\end{Corollary}

\begin{Theorem}
  A pair of linear maps $(\delta_+,\delta_-)$ on a Jordan *-algebra is a Jordan triple
  derivation pair if, and only if,
  \begin{enumerate}[(i)]
  \item $\delta_\pm(a\j b) = \delta_\pm a\j b + a\j \delta_\pm b + \tp
    a{\delta_\mp1}b$ for all $a$ and $b$, and
  \item $\delta_\pm(b^*) = 2\delta_\pm1\j b^* +(\delta_\mp b)^*$ for all $b$.
  \end{enumerate}

\end{Theorem}

\section{Local triple derivations}

In~\cite{MR1051316}, Kadison defines a local derivation to be a
linear map which at each point $a$ takes the same value as some
derivation.  This would appear to be a significant generalisation of
derivation (an example, due to C.~Jensen, of a local derivation which
is not a derivation is given) but Kadison
proceeds to show that, in the setting of  von Neumann
algebras at least, the definition is void: 
\begin{Theorem}[{\cite{MR1051316}}, Thm A]
  Every continuous local derivation on a von Neumann algebra is a derivation.
\end{Theorem}
The proof is rather lengthy and ingenious.  The result has the
following immediate corollary, notable enough to be labelled as a
theorem.
\begin{Theorem}[{\cite{MR1051316}}, Thm B]\label{thm:B}
  If $\delta$ is a norm continuous linear mapping of a von Neumann
  algebra into itself such that for each $a\in  A$ there exists $x_a
  \in A$ with $\delta(a) = [a, x_a]$ then there exists $x\in A$ with
  $\delta = [\cdot, x]$.
\end{Theorem}
Our aim
in this section 
is to introduce a similar definition of \emph{local triple derivation}
and to prove the analogous result that every continuous local triple derivation
on a JBW*-triple is a triple derivation.  An analogue of \ref{thm:B}
will ensue.  It is only fair to point out
that, as we have seen earlier, Jordan triple derivations are not
generalisations of algebra derivations, and so our result runs
parallel to Kadison's rather than being a generalisation of it.  
Even so, it is perhaps surprising that our proof has little in common
the binary case and is substantially more compact.  We attribute this
to the elegant symmetry 
of the Jordan setting rather than any fresh ingenuity.

\begin{Definition}\label{def:ltd}
  A \textsc{local triple derivation} on a Jordan triple $Z$
  is a linear map $\delta:Z\to Z$ such
  that, for every $x\in Z$, there exists a derivation $\delta_x$ with
  $\delta_x(x)=\delta(x)$. 
\end{Definition}

Of course, every triple derivation is a local triple derivation and
a natural question is whether there exist local derivations which are
not derivations.  The following example, a variation of one
attributed to C.U.~Jensen in \cite{MR1051316}, shows that such maps do exist
in a purely algebraic setting.

\begin{Example}
  Consider the *-algebra $\mathbb{C}(x)$ of rational functions in the variable $x$
  over $\mathbb{C}$ and its *-subalgebra of polynomials
  $\mathbb{C}[x]$.  This algebra, as pointed out in~\cite{MR1051316},
  provides an example of a local derivation which is not a
  derivation.  The derivations of $\mathbb{C}(x)$ take the form
  $d_g:=f\mapsto gf'$ for $g\in \mathbb{C}(x)$.  Note that since the algebra is
  commutative, derivations and Jordan derivations agree.  We are
  interested in triple derivations and we make the following remarks:
  \begin{enumerate}[(i)]
  \item Either by direct calculation, or appealing to
    Theorem~\ref{thm:char}, one finds that a linear map $\delta$ on
    $\mathbb{C}(x)$ is a triple derivation if 
    $\delta=\delta_{u,v}$ where $\delta_{u,v}f= uf'+ivf$ for
    $f\in\mathbb{C}(x)$ and $u,v$ are self-adjoint
    elements of $\mathbb{C}(x)$.  
    \item All triple derivations of $\mathbb{C}(x)$ are of the form
      $\delta_{u,v}$.  To see this, let $\delta$ be a triple
      derivation.  By Theorem~\ref{thm:char}, $(\delta 1)^*=-\delta1$
      and thus $\delta 1 =iv$ for some $v=v^*$.  Also
      by~\ref{thm:char} since $x=x^*$, $\delta x = 2(\delta1) x +(\delta
      x)^*$.  Thus, if we let $u=\delta x - x\delta 1$ then $u^*=u$.
      Now, one more use of~\ref{thm:char}(a) shows that 
      \begin{align*}
        \delta(x^2)&= 2x\delta x + \tp x{\delta1}x \\
                  &= 2x\delta x-x^2\delta1\\
                  &=2x(\delta x -x\delta1) + x^2\delta 1 \\
                  &= u(x^2)' + iv(x^2).
                \end{align*}
     It is but a short step to showing $\delta(x^n)= u(x^n)' +
     iv(x^n)$ and hence $\delta(p)=up'+ivp$ for all $p\in
     \mathbb{C}[x]$.  The extension to all $f\in \mathbb{C}(x)$
     follows.
    \item The \emph{local}
    triple derivations are the linear maps which take $1$ to $iv$
    where $v=v^*\in \mathbb{C}(x)$.  Indeed, if $\alpha$ is a local
    triple derivation then $\alpha(1)=\delta_{u,v}(1)$ for some
    self-adjoint $u$ and $v$.  But $\delta_{u,v}(1)=iv$.  Conversely,
    if $\alpha$ is linear with $\alpha(1)=iv$ then $\alpha$ agrees
    with the derivation $\delta_{0,v}$ at any constant.  Solving the
    functional equation $uf' + ivf = \alpha(f)$ for
    self-adjoint functions $u$ and $v$ leads to 
\[ u = \frac{\psi_1\phi_1 + +\psi_2\phi_2}{\phi_1'\phi_1 +
  \phi_2'\phi_2}, \qquad v= \frac{\psi_2\phi_1'-\psi_1\phi_2}{\phi_1'\phi_1 +
  \phi_2'\phi_2}\] where $\alpha(f)=\psi_1+i\psi_2$
and $f=\phi_1+i\phi_2$ and these solutions exist in $\mathbb{C}(x)$ as
long as $(\phi_1^2 + \phi_2^2)'$ is not zero, that is, when $ff^*$
(and consequently $f$) is not constant.  Thus $\alpha(f)=d_f(f)$ for
some derivation $d_f$ whether $f$ is constant or not, and so $\alpha$ is
a local derivation.
\item The linear map $f\mapsto i(xf)'$, which can be written here as
  $\delta_{ix,1}$, is now seen to be a local derivation ($1\mapsto
  i1$) but not a derivation ($ix \ne (ix)^*$). 
  \end{enumerate}
\end{Example}

Having seen that there do exist local triple derivations which are not
derivations, let us proceed to show that no such examples exist in any
JBW*-triple.

\begin{Lemma}\label{lem:54}
  Let $e$ and $f$ be orthogonal tripotents and $\delta$ a local
  derivation on a Jordan triple.  Then
  \[ \delta\tp fef = 2\tp{\delta(f)}ef + \tp f{\delta e}f.\]
\end{Lemma}

\begin{Proof}
  As $\D ef =0$ one need only show that $\tp f{\delta e}f = 0$.
For this, choose a derivation $\delta_e$ with $\delta e= \delta_e e$.
Since $\delta_e\tp fef = 2\tp{\delta_e(f)}ef + \tp f{\delta_e e}f$, we
may again ignore zero products to conclude $0= \tp f{\delta_e e}f$ which
gives the result.
\end{Proof}

We seek a similar result where the triple product is of the form $\tp
ffe$.

\begin{Lemma}
  Let $e$ and $f$ be orthogonal tripotents in a Jordan triple $Z$ and $\delta$ a local
  derivation on $Z$.  Then $\tp e{\delta e}f +\tp ee{\delta f}=0.$
\end{Lemma}

\begin{Proof}
  For any tripotent $e$ and local derivation $\delta$, we have
  $\delta(e)=\delta_e(e)=\delta_e\tp eee=2\tp{\delta_ee}ee+\tp
  e{\delta_ee}e= 2\tp {\delta e}ee + \tp e{\delta e}e$.  Now, if $e$
  and $f$ are orthogonal tripotents then $e+f$ is also a tripotent and
  so
\[ \delta(e+f)= 2\tp{\delta(e+f)}{e+f}{e+f} +\tp {e+f}{\delta(e+f)}{e+f}.\]
This leads to \[ 2\tp {\delta e}ff + 2\tp {\delta f}ee + 2\tp e{\delta
  e}f + \tp e{\delta f}e + 2\tp e{\delta f}f + \tp
f{\delta e}f  =0\]
and, after replacing $e$ by $-e$ and summing, we conclude \[ 2\tp{\delta f}ee + 2\tp
e{\delta e}f + \tp e{\delta f}e =0.\]
From the proof of Lemma~\ref{lem:54}, $\tp e{\delta f}e=0$ and thus \[
\tp {\delta f}ee + \tp e{\delta e} f =0\] as asserted.
\end{Proof}

\begin{Corollary}\label{cor:58}
  Let $e$ and $f$ be orthogonal tripotents on a Jordan triple $Z$ and $\delta$ a local
  derivation on $Z$.  Then \[ \delta \tp eef = \tp{\delta e}ef + \tp e{\delta e}f
  + \tp ee{\delta f}.\]
\end{Corollary}

\begin{Lemma}\label{lem:59}
  Let $e$, $f$ and $g$ be mutually orthogonal tripotents on a Jordan
  triple $Z$ and $\delta$
  a local derivation on $Z$.  Then  \[ \delta\tp efg = \tp{\delta e}fg + \tp
  e{\delta f}g + \tp ef{\delta g}.\]
\end{Lemma}

\begin{Proof}
  Excluding the zero products, one must only show that $\tp e{\delta
    f}g=0$.  For this, one need only choose a derivation $\delta_f$
  such that $\delta_f(f)=\delta f$ and remark that $\tp e{\delta_f
    f}g=0$ since the desired identity holds for a derivation.
\end{Proof}

At this point, we have effectively considered a number of different
cases which are covered by the following proposition.

\begin{Proposition}
  Let $\delta$ be a local triple derivation on a Jordan triple $Z$, $\Lambda$ a family of
  orthogonal tripotents and $e$, $f$ and 
  $g$ elements of $\Lambda$ (not necessarily distinct).  Then 
  \[
    \delta\tp efg = \tp {\delta e}fg + \tp e{\delta f}g + \tp
    ef{\delta g}.
  \]
\end{Proposition}

\begin{Proof}
  If $e=f=g$ then the result follows on replacing $\delta$ by a
  derivation $\delta_e$ with $\delta_e(e)=\delta e$.  If $e$, $f$ and
  $g$ are all distinct then the result holds by Lemma~\ref{lem:59}.
  If, on the other hand, there are just
  two distinct tripotents, then our triple product is
  either of the form $\tp efe$ or $\tp eef$ and the conclusion is reached
  on appeal to \ref{lem:54} or \ref{cor:58} as appropriate.
\end{Proof}

This can be further extended via the linearity of our local
derivation.
\begin{Corollary}\label{cor:alg}
  Let $\Lambda$ be a family of orthogonal tripotents and $x$, $y$ and
  $z$ elements of ${\mathrm{span}\ \Lambda}$.  Then 
  \[
    \delta\tp xyz = \tp {\delta x}yz + \tp x{\delta y}z + \tp
    xy{\delta z}.
  \]
\end{Corollary}

We can now present our main result.

\begin{Theorem}\label{thm:main}
  Let $\delta$ be a continuous local triple derivation on a JBW*-triple.  Then $\delta$
  is a triple derivation.
\end{Theorem}

\begin{Proof}
  Fix an element $x$ and choose $\epsilon \in (0,1)$.  Then, by Proposition~\ref{prop:dens}, we can find mutually
  orthogonal tripotents $\{e_k:k=0,\ldots,n\}$, a linear combination
  of which, say $\xi = \summ k0n \alpha_k e_k$,
  has the property that
  $\norm{x-\xi}<\epsilon$.    Corollary~\ref{cor:alg} asserts
  that $\delta\tp\xi\xi\xi = 2\tp{\delta\xi}\xi\xi +
  \tp\xi{\delta\xi}\xi$.  It follows from $\norm{\tp
    yzw}\le \norm y\norm z\norm w$ that
  \( \norm{\tp xxx -\tp
  {\xi}{\xi}{\xi}}\le C_1 \norm{x-\xi}< C_1\epsilon \)
and thus by continuity of $\delta$ that
\( \norm{\delta(\tp xxx) - \delta(\tp
  {\xi}{\xi}{\xi})}\le C_2\epsilon\).
In a similar vein, $\norm{\tp{\delta x}xx - \tp{\delta
  \xi}{\xi}{\xi}}\le C_3\epsilon$ and $\norm{\tp x{\delta x}x - \tp
{\xi}{\delta \xi}{\xi}}\le C_4\epsilon$ where each $C_j$ is positive
and independent of $\epsilon$. 
Let us now estimate as follows:
\begin{align*}
  \norm{\delta\tp xxx - (2\tp {\delta x}xx + \tp x{\delta x}x)} &\le
     \norm{\delta\tp xxx - \delta\tp \xi\xi\xi} \\
     &\qquad + \norm{\delta\tp\xi\xi\xi - (2\tp {\delta \xi}\xi\xi +
       \tp \xi{\delta \xi}\xi)} \\
     &\qquad + \norm{2\tp {\delta \xi}\xi\xi - 2\tp{\delta x}xx} \\
     &\qquad +\norm{\tp \xi{\delta \xi}\xi - \tp x{\delta x}x}\\
     &\le C_2\epsilon + 0 + 2C_3\epsilon + C_4\epsilon.
\end{align*}
As $\epsilon$ is arbitrarily small, we can conclude that
\begin{equation}
  \label{eq:hom}
  \delta\tp xxx = 2\tp{\delta x}xx + \tp x{\delta x}x.
\end{equation}
The final step is a polarisation exercise.  Replacing $x$ first by
$x+y$, then by $x-y$ in \eqref{eq:hom} and summing leads to
\begin{align}
  \delta\tp xyy +\delta \tp yxy +\delta \tp yyx &= 2(\tp {\delta x}yy + \tp
  {\delta y}xy + \tp {\delta y}yx) \notag\\ &\qquad\qquad+ \tp x{\delta y}y +\tp y{\delta
    x}y + \tp y{\delta y}x.  \label{eq:quad}
\end{align}
In \eqref{eq:quad}, replace $x$ by $ix$ and compare with
\eqref{eq:quad} multiplied by $i$ to gain
\begin{equation}
  \label{eq:quad1}
  \delta\tp yxy = 2\tp{\delta y}xy + \tp y{\delta x}y.
\end{equation}
Finally, replace $y$ by $y+z$ in \eqref{eq:quad1} to conclude that
\( \delta\tp yxz = \tp {\delta y}xz + \tp y{\delta x}z + \tp yx{\delta
  z}\)
and $\delta$ is a derivation.
\end{Proof}

Notice that continuity of the local derivation was used in quite a
weak sense in this proof.  If each element of the JBW*-triple is
represented by a \textsl{finite} linear combination of orthogonal
tripotents, then the conclusion remains true without appealing to the
continuity of the local derivation.  This is the case in any finite
rank JBW*-triple and so we have the following.

\begin{Theorem}
  Let $\delta$ be a local derivation on a finite rank JBW*-triple.
  Then $\delta$ is a derivation (and hence is continuous).
\end{Theorem}

As an application of Theorem~\ref{thm:main}, we present the following.
\begin{Corollary}
  Suppose $r$ is a bounded linear map on a von Neumann algebra $A$ and
  for every $x\in A$ there exists $p_x\ge 0$ such that $r(x)=xp_x+p_x
  x$.  Then there exists $p\ge 0$ such that for every $x\in A$,
  \[ r(x^2)=xr(x)+r(x)x - xpx.\]
\end{Corollary}

\begin{Proof}
  Writing $p_x=h_x^2$ for $h_x=h_x^*$ we see that $ir(x)= i(xh_x^*h_x
  + h_xh_x^*x)=i\D{h_x}{h_x} x$.  That is, $ir$ agrees with a triple
  derivation at every point and so is a local triple derivation.
  Therefore, by Theorem~\ref{thm:main}, $ir$ is a triple derivation
  and, in
  particular,
  \begin{align*}
    ir(x^2)= ir\tp x1x &= 2\tp{irx}1x + \tp x{ir(1)}x \\
                       &= 2i(rx \j x) + x(ip)^*x\\
                       &= i(x.rx + rx.x)-ixpx
  \end{align*}
  where $p=r(1)\ge 0$.  
\end{Proof}

\begin{Corollary}
  Suppose $s$ is a bounded linear map on a von Neumann algebra $A$ and
  for every $x\in A$ there exists $p_x\ge 0$ such that either
  $(s(x),s(x^*))=(xp_x,x^*p_x)$ or $(s(x),s(x^*))=(p_xx,p_xx^*)$.  Then there exists $p\ge 0$ such that
  for every $x\in A$,
  \[ s(x^2)=xs(x)+s(x)x - xpx.\]
\end{Corollary}

\begin{Proof}
  Defining $s^*(x)=(s(x^*))^*$ we have $s$ is a bounded linear map and
  $(s+s^*)(x)=p_x x+xp_x$.  Therefore, we can apply the previous corollary to
  $s+s^*$ to gain
  \[ (s+s^*)(x^2)= x(s+s^*)(x)+(s+s^*)(x)x - xpx.\]
  Also, $(s-s^*)(x)= \pm(xp_x-p_xx)$ which means that, at $x$, $s-s^*$ agrees
  with an inner derivation.  By
  Kadison's original result, $s-s^*$ is a
  derivation.  Thus
  \[ (s-s^*)(x^2) = x(s-s^*)(x) + (s-s^*)(x)x.\]
  Summing we see $s(x^2)= xs(x)+s(x)x - xpx$ as required.
\end{Proof}

We can go further in this direction.  Notational changes are enough to extend
Theorem~\ref{thm:main} to the following variant for derivation pairs.
We call a pair of linear maps $(\delta^+,\delta^-)$ a \emph{local
  derivation pair} if for every $x$ the  re exists a derivation pair
$(d_x^+, d_x^-)$ such that $\delta^+(x)=d^+(x)$ and
$\delta^-(x)=d^-(x)$.

\begin{Theorem}
  Every local derivation pair on a JBW*-triple is a derivation pair.
\end{Theorem}

\section{Concluding remarks and open questions}

B.E. Johnson \cite{MR1783788} provided the following strong extension of
Kadison's theorem:
\begin{Theorem}\label{thm:Johnson}
  Every local derivation on a C*-algebra is a derivation.
\end{Theorem}
Notice that, apart from widening the class of algebras dealt with,
Johnson's result drops the requirement of continuity of the local
derivation.  In particular, the automatic continuity of local C*-algebra
derivations follows from the automatic continuity of C*-algebra
derivations.  This raises natural conjectures in the triple setting.

\begin{Conjecture}
  \begin{enumerate}[({C}1)]
      \item \label{conj:1} A local triple derivation on a JB*-triple is a derivation.
   \item \label{conj:2} A local triple derivation on a JB*-triple is continuous. 
 \item \label{conj:3} A continuous local triple derivation on a JB*-triple is a
    derivation.
  \end{enumerate}
\end{Conjecture}

\newcommand\cref[1]{(C\ref{#1})}

Clearly \cref{conj:1} (or rather, proof thereof) implies \cref{conj:3}
and, by the automatic continuity of derivations \cite{MR1067482},
\cref{conj:1} also implies \cref{conj:2}.  Conversely, \cref{conj:2}
and \cref{conj:3} together imply \cref{conj:1}.  Johnson's proof
of \ref{thm:Johnson} relies on use of the multiplier algebra of a C*-algebra.

The reader will bear in mind that the results of Kadison and Johnson
were proven for \textsl{module-valued} derivations on von Neumann and
C*-algebras respectively, while in this paper we have restricted to
consideration of triple derivations of the Jordan triple into
itself. Recently, Peralta and Russo \cite{PeraltaRusso_derivations} have initiated a study of
module-valued Jordan triple derivations.  In particular, they answer
the question of when a (module-valued) Jordan triple derivation is
automatically continuous.

\medskip
\noindent \textbf{Remark.} A triple automorphism on a Jordan triple is
a bijective linear map $T$ satisfying $T\tp xyz=\tp{Tx}{Ty}{Tz}$.  The
ensuing definition of local triple automorphism, analogous to
Definition~\ref{def:ltd}, is clear: a linear map $T$ is a local triple
automorphism if, for every $x$, there exists a triple automorphism
$T_x$ such that $T(x)=T_x(x)$.  For any JB*-triple, the set of
triple automorphisms and the set of surjective linear isometries coincide.  Thus,
if $T$ is a local triple automorphism then it is an isometry since, for each $x$, $\norm
{Tx}=\norm {T_x(x)}= \norm x$.  In
particular, we see immediately  that (a) any local triple automorphism is continuous
and (b) any surjective local triple automorphism is a surjective
linear isometry and hence a triple automorphism.

\bigskip

\noindent \textbf{Acknowledgement.} The author thanks R.V.~H\"ugli for
helpful comments.


\bibliographystyle{acm}
\bibliography{/home/mackey/research/bib}

\def\cprime{$'$}
\begin{thebibliography}{10}

\bibitem{BT}
{\sc Barton, T., and Timoney, R.~M.}
\newblock Weak{$^*$}-continuity of {J}ordan triple products and its
  applications.
\newblock {\em Math. Scand. 59\/} (1986), 177--191.

\bibitem{MR1067482}
{\sc Barton, T.~J., and Friedman, Y.}
\newblock Bounded derivations of {${\rm JB}^*$}-triples.
\newblock {\em Quart. J. Math. Oxford Ser. (2) 41}, 163 (1990), 255--268.

\bibitem{Cartan}
{\sc Cartan, E.}
\newblock Sur les domaines born\'es homog\`enes de l'espace de $n$ variables
  complexes.
\newblock {\em Abh. Math. Semin. Univ. Hamburg 11\/} (1935), 116--162.

\bibitem{MR2190348}
{\sc Chu, C.-H., and Mackey, M.}
\newblock Isometries between {${\rm JB}^*$}-triples.
\newblock {\em Math. Z. 251}, 3 (2005), 615--633.

\bibitem{Dineen_CHVF}
{\sc Dineen, S.}
\newblock Complete holomorphic vector fields on the second dual of a {B}anach
  space.
\newblock {\em Math. Scand. 59\/} (1986), 131--142.

\bibitem{FR_GN}
{\sc Friedman, Y., and Russo, B.}
\newblock The {G}elfand {N}aimark theorem for {JB*-triple}s.
\newblock {\em Duke Math. J. 53}, 1 (1986), 139--148.

\bibitem{HOS}
{\sc Hanche-Olsen, H., and St\"ormer, E.}
\newblock {\em Jordan operator algebras}.
\newblock Pitman, 1984.

\bibitem{MR1875137}
{\sc Ho, T., Martinez-Moreno, J., Peralta, A.~M., and Russo, B.}
\newblock Derivations on real and complex {JB} {$^\ast$}-triples.
\newblock {\em J. London Math. Soc. (2) 65}, 1 (2002), 85--102.

\bibitem{MR929400}
{\sc Horn, G.}
\newblock Characterization of the predual and ideal structure of a {${\rm
  JBW}^*$}-triple.
\newblock {\em Math. Scand. 61}, 1 (1987), 117--133.

\bibitem{MR2491603}
{\sc H{\"u}gli, R.~V., and Mackey, M.}
\newblock Transitivity of inner automorphisms in infinite dimensional {C}artan
  factors.
\newblock {\em Math. Z. 262}, 1 (2009), 125--141.

\bibitem{MR1997703}
{\sc Isidro, J.~M.}
\newblock Holomorphic automorphisms of the unit balls of {H}ilbert
  {$C^*$}-modules.
\newblock {\em Glasg. Math. J. 45}, 2 (2003), 249--262.

\bibitem{MR1783788}
{\sc Johnson, B.~E.}
\newblock Local derivations on {$C^*$}-algebras are derivations.
\newblock {\em Trans. Amer. Math. Soc. 353}, 1 (2001), 313--325.

\bibitem{MR1051316}
{\sc Kadison, R.~V.}
\newblock Local derivations.
\newblock {\em J. Algebra 130}, 2 (1990), 494--509.

\bibitem{Kaup_RMT}
{\sc Kaup, W.}
\newblock A {R}iemann mapping theorem for bounded symmetric domains in complex
  {B}anach spaces.
\newblock {\em Math. Z. 138\/} (1983), 503--529.

\bibitem{KaupUpmeier_Siegel}
{\sc Kaup, W., and Upmeier, H.}
\newblock Jordan algebras and symmetric {S}iegel domains in {B}anach spaces.
\newblock {\em Math. Z. 157\/} (1977), 179--200.

\bibitem{Mackey_homotopes}
{\sc Mackey, M.}
\newblock Homotopes of {JB*}-triples.
\newblock {\em Asian-European Journal of Mathematics 2}, 3 (2009), 465--475.

\bibitem{McCrimmon_taste}
{\sc McCrimmon, K.}
\newblock {\em A taste of {J}ordan algebras}.
\newblock Universitext. Springer-Verlag, New York, 2004.

\bibitem{PeraltaRusso_derivations}
{\sc Peralta, A.~M., and Russo, B.}
\newblock Automatic continuity of derivations on {C*}-algebras and
  {JB*}-triples.
\newblock {\em Preprint\/} (2012).

\end{thebibliography}

\end{document}